\documentclass[a4paper,12pt]{article}
\usepackage[text={165mm,245mm}]{geometry}
\usepackage{graphicx} 
\usepackage[utf8]{inputenc}
\usepackage{authblk}
\usepackage{subfigure}
\usepackage{float}
\usepackage{amsmath}
\usepackage{amsthm}
\usepackage{amssymb}
\usepackage{mathtools}
\usepackage{color}

\usepackage[hyperfootnotes=true,colorlinks=true]{hyperref}
\hypersetup{
    colorlinks=true,
    linkcolor=blue,
    citecolor=magenta,      
    urlcolor=black
}

\title{The M\"obius--Kantor graph is a faithful \\ unit-distance graph}

\author[1,2,3]{Nino Ba\v{s}i\'{c}}
\author[4]{G\'{a}bor G\'{e}vay}
\author[1,2,3]{Toma\v{z} Pisanski}

\affil[1]{FAMNIT, University of Primorska, Koper, Slovenia}
\affil[2]{IAM, University of Primorska, Koper, Slovenia}
\affil[3]{Institute of Mathematics, Physics and Mechanics, Ljubljana, Slovenia}
\affil[4]{Bolyai Institute, University of Szeged, Szeged, Hungary}

\date{\today}

\begin{document}

\maketitle

\begin{center}
\textbf{\emph{To Tom Tucker on the happy occasion of his 80th birthday.}}
\end{center}

\begin{abstract}
In this paper, it has been shown that the generalized Petersen graph $\mathrm {GP}(8,3)$, also known as the 
M\"obius--Kantor graph, admits a faithful unit-distance representation in the plane. 


\bigskip
\noindent \textbf{Math.\ Subj.\ Class.\ (2020):} 
05C10,  %  Planar graphs; geometric and topological aspects of graph theory
05C62,  % Graph representations (geometric and intersection representations, etc.) 
05E18,  % Group actions on combinatorial structures
20B25.  % Finite automorphism groups of algebraic, geometric, or combinatorial structures
\end{abstract}

The generalized Petersen graph $\mathrm {GP}(8,3)$, 
also known as the M\"obius--Kantor graph, appears in
many mathematical themes \cite{MaPi2000}. Its name
comes from the fact that it is the Levi graph of the 
M\"obius--Kantor $(8_3)$ configuration
\cite{PiSe2013}. Its automorphism group is of order 96
and is the only group of genus two, as shown by Thomas W.\ Tucker 
\cite{Tu1984}. 

\begin{figure}[!htbp]
\centering
    \subfigure[]{ \hskip -8pt
      \includegraphics[width=0.47\linewidth]{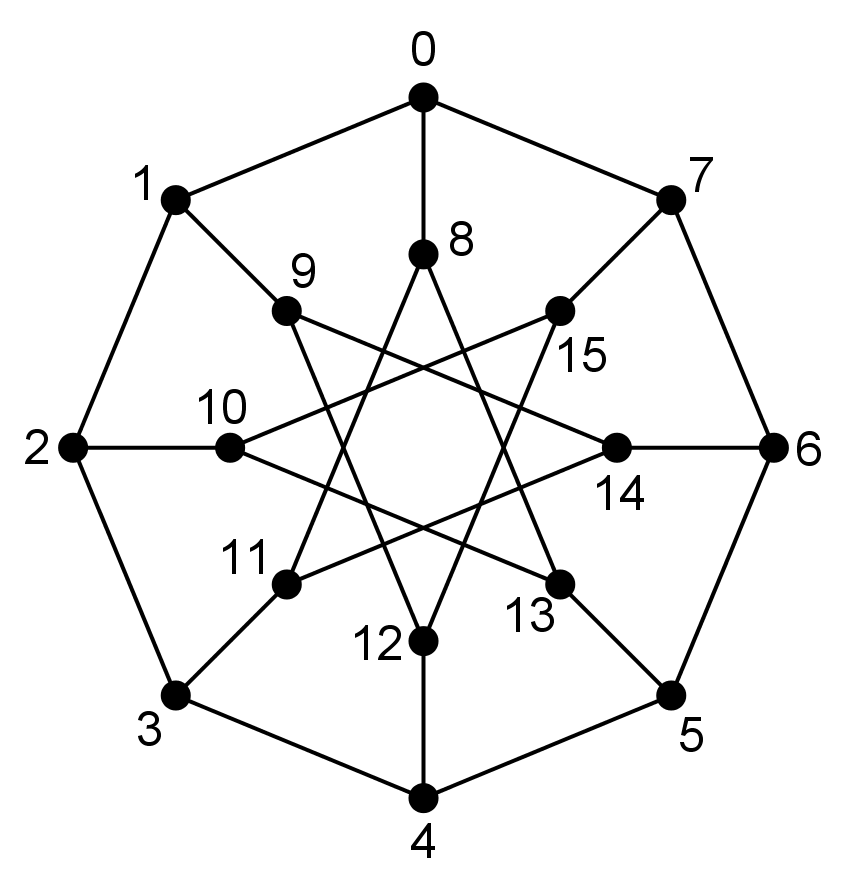}} \hskip 15 pt
    \subfigure[]{ \hskip -10pt
      \includegraphics[width=0.445\linewidth]{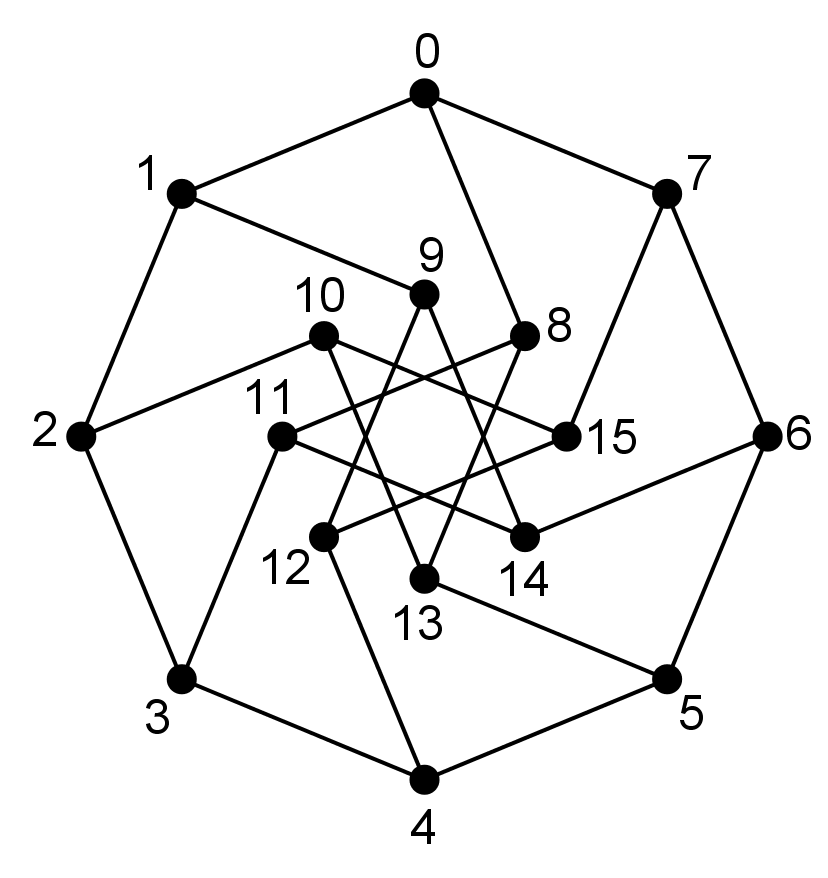}}     
\caption{The M\"obius--Kantor graph. (a) Standard drawing with 
vertex labels as used in Table~\ref{eq:theMainSystem}. 
(b) Unit-distance representation that is not faithful. For each 
vertex there is a non-adjacent vertex at a distance 1 from it; 
for example, the distance between vertices 0 and 10 is 1.}
  \label{fig:MK_usual}
\end{figure}

It is well known that the M\"obius--Kantor graph, as
any other generalized Petersen graph, is a unit-distance graph
\cite{ZiHoPi2012}. The drawing depicted in 
Figure~\ref{fig:MK_usual}(b) is usually taken as a prime example of a 
unit-distance representation that is not \textit{faithful} or 
\emph{strong}. In practice, this means that all edges have length 1, 
however, there exist pairs of non-adjacent vertices that are placed at 
a distance 1. 

A formal definition, taken from \cite{NoKu2014}, is as 
follows. A \emph{faithful unit-distance graph} in $\mathbb R^d$ 
is a graph whose set of vertices is a finite subset of the 
$d$-dimensional Euclidean space, where two vertices are adjacent 
if and only if the Euclidean distance between them is exactly 1. 
A \emph{unit-distance graph} in $\mathbb R^d$ is any subgraph of 
such a graph.

Note that one has to be precise when addressing possible 
degeneracies \cite{HoPi2012,PiZi2009}. That is, even in the 
case of faithful unit-distance graphs, a vertex may be mapped 
to the interior of an interval representing an edge to which
it does not belong. In general, there exist unit-distance graphs
with a partial overlap of certain represented edges. 

\begin{figure}[!htbp]
    \centering
    \includegraphics[width=0.575\linewidth]
    {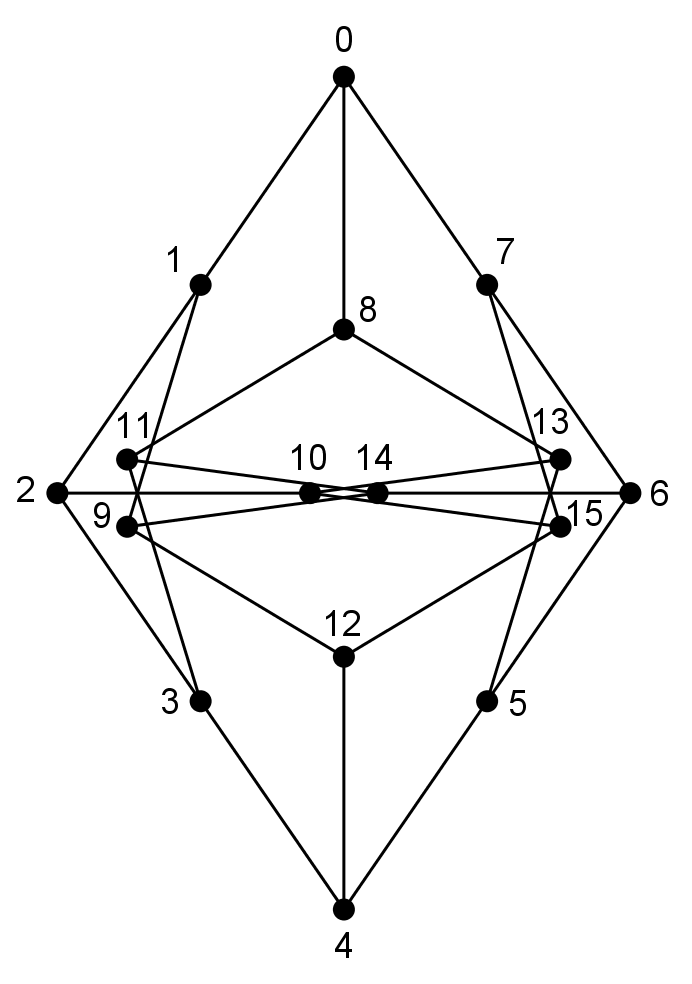}
    \caption{Faithful unit-distance representation of the 
    M\"obius--Kantor graph. The vertex labels refer to our 
    Table~\ref{tbl:coordinatesOfVertices}; also cf.\ Figure 
    \ref{fig:MK_usual}(a).}
    \label{fig:faithful}
\end{figure}

Recently, we discovered
the representation of $\mathrm{GP}(8,3)$ that is shown
in Figure~\ref{fig:faithful}. It can be proven that this
representation is a unit-distance representation with a 
rhombus-shaped outermost 8-cycle; moreover, it can easily 
be checked that the distances between its non-adjacent 
vertices are not close to 1. Hence, the representation 
depicted in Figure~\ref{fig:faithful} is in fact a faithful
unit-distance representation. It turns out that there 
exist infinitely many such representations. The context of 
this finding will be published elsewhere.

Note that the novel representation has a dihedral symmetry 
of order 4 and is polycirculant with respect to the group
${\mathbb Z}_2$; compare \cite{BeGePi2025}.

If we insist on that the symmetry is dihedral
${\mathbb D}_2$ and require that the outer shape be a 
rhombus, choose the following coordinates for three vertex 
representations:
\begin{align*}
0 &\colon (0,k), \\
6 &\colon (h,0), \\
13 &\colon (p,q),
\end{align*}
which then define all other vertex coordinates.
The values $h, k, p$ and $q$ are determined by the following system of quadratic 
equations, each of them obtained by using Pythagoras' 
theorem:
\begin{equation}
\begin{aligned}
h^2 + k^2 & = 2^2,\\
p ^2 +(q-k +1)^2  & = 1, \\
q ^2 +(p+h -1)^2  & = 1, \\
(p - h/2)^2 + (q + k/2)^2 & = 1. \\
%eq1 = p^2 + (q - k + 1)^2 == 1
%eq2 = (p + h - 1)^2 + q^2 == 1
%eq3 = (2*p - h)^2 + (2*q + k)^2 == 4
%eq4 = k^2 + h^2 == 4
\end{aligned}
\label{eq:theMainSystem}
\end{equation}
The full set of coordinates is given in Table~\ref{tbl:coordinatesOfVertices}.

\begin{table}[!htbp]
\centering
\begin{tabular}{|c|l|}
\hline
Vertex & Coordinates \\
\hline
\hline
0 & $(0, k)$ \\
1 & $(-h/2, k/2)$ \\
2 & $(-h, 0)$ \\
3 & $(-h/2, -k/2)$ \\ 
4 & $(0, -k)$ \\
5 & $(h/2, -k/2)$ \\ 
6 & $(h, 0)$ \\
7 & $(h/2, k/2)$ \\ 
8 & $(0, k-1)$ \\
9 & $(-p, -q)$ \\
10 & $(1-h, 0)$ \\
11 & $(-p, q)$ \\
12 & $(0, 1-k)$ \\
13 & $(p, q)$ \\
14 & $(h-1, 0)$ \\ 
15 & $(p, -q)$ \\
\hline
\end{tabular}
\caption{The coordinates of the 16 vertices of the graph $\mathrm{GP}(8, 3)$, where $h, k, p$ and $q$ satisfy
the system of equations~\eqref{eq:theMainSystem}.}
\label{tbl:coordinatesOfVertices}
\end{table}

The system has two non-degenerate real solutions. One solution is:
\begin{align*}
h & \approx 1.133693, \\
k & \approx 1.647647, \\
p & \approx 0.857420, \\
q & \approx 0.133029,
\end{align*}
with the corresponding drawing presented in Figure~\ref{fig:faithful}. 
The other is just a reflection in the line $y = x$.

\bigskip
As we mentioned above, it is easy to verify that all
edges of ${\mathrm{GP}(8,3)}$ are indeed of length 1 and
that no non-adjacent vertices are at a distance 1 from
each other.

We conclude our paper by briefly mentioning a novel consequence
of our result in the context of configurations, which follows 
from the construction invented in \cite{GePi2014}.

\begin{figure}[!htbp]
\centering
    \subfigure{ \hskip -10pt
      \includegraphics[width=0.415\linewidth]{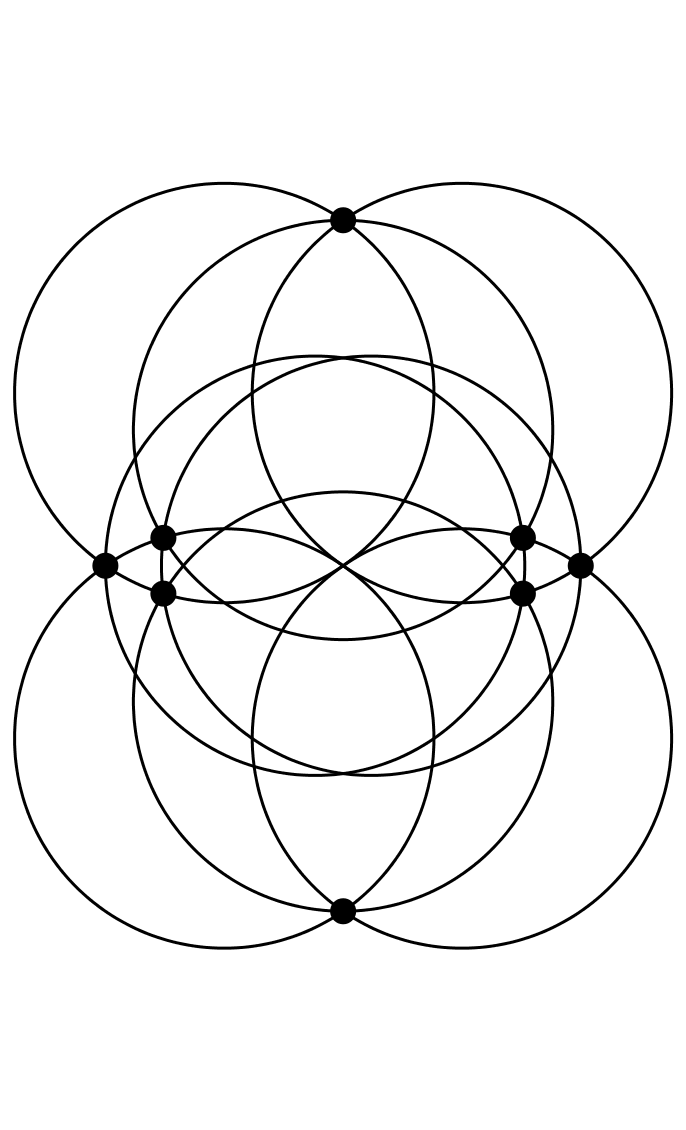}} \hskip 5 pt
    \subfigure{ \hskip -6pt
      \includegraphics[width=0.545\linewidth]{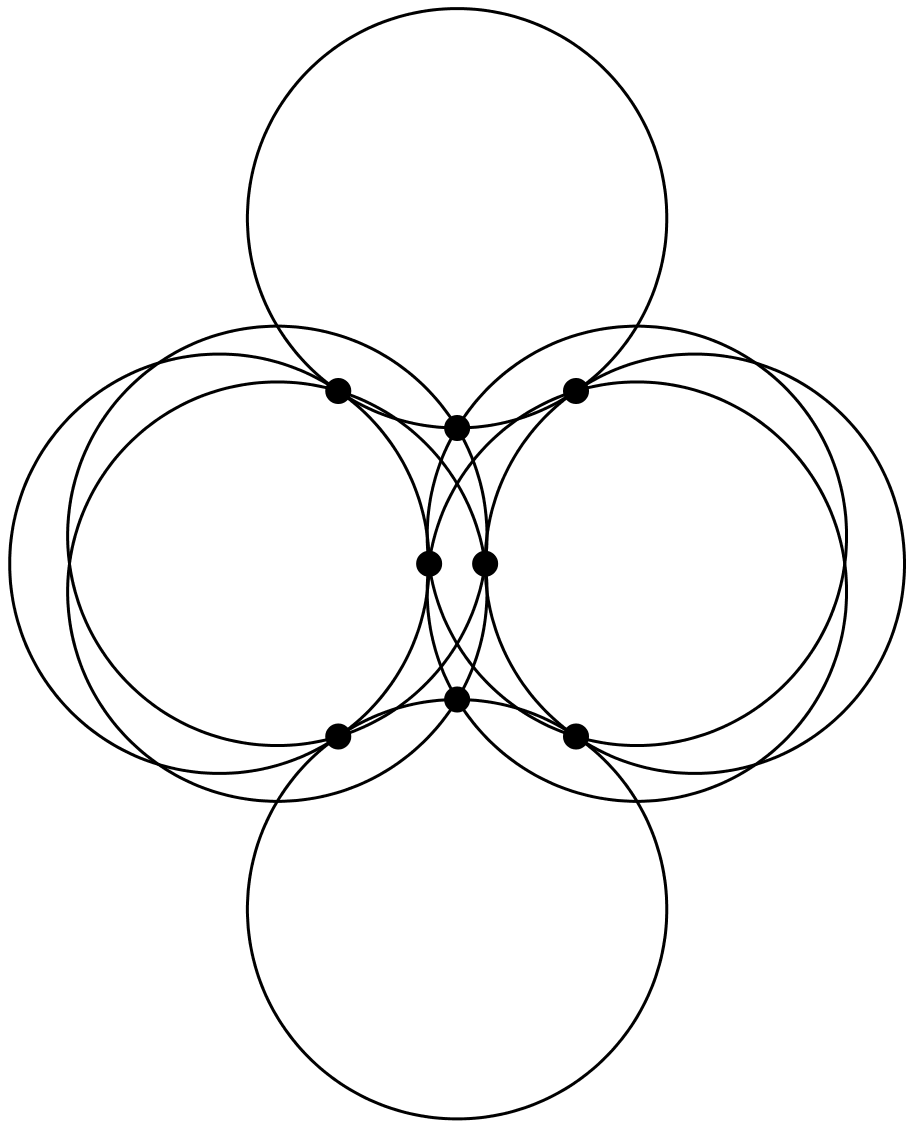}}
\caption{Isometric point-circle realizations of the 
M\"obius--Kantor configuration derived from the faithful unit-distance 
representation of the M\"obius--Kantor graph in 
Figure~\ref{fig:faithful}. They are dual to each other.}
\label{fig:isom}
\end{figure}

It is well known 
that the M\"obius--Kantor configuration cannot be geometrically 
realized with points and (straight) lines in the real
Euclidean plane~\cite{PiSe2013}. However, a faithful unit-distance
representation of the M\"obius--Kantor graph provides the 
possibility to realize it with points and circles, in the 
following way. Since the graph is bipartite and unit-distance, any 
of the bipartition classes can be taken as the set of centres of 
unit circles, while the other class serves as the set of 
configuration points. Furthermore, since the graph 
is faithful, this guarantees that no false incidences occur in the 
configuration. By exchanging the role of the two bipartition 
classes one obtains a second configuration such that the two 
configurations are dual to each other (see 
Figure~\ref{fig:isom}). 
As the radius for all the circles is of unit length, we 
call such configurations \emph{isometric} point-circle 
configurations~\cite{GePi2014}.

\section*{Acknowledgements}

Nino Ba\v{s}i\'{c} is supported in part by the Slovenian Research Agency (research program P1-0294 and research
project J5-4596).
Toma\v{z} Pisanski is supported in part by the Slovenian Research Agency (research program P1-0294 and research
projects J1-4351, J5-4596, and BI-HR/23-24-012).

% \bibliography{MKFUD}
% \bibliographystyle{amcjoucc} 

\end{document}